\documentclass{ifacconf}

\usepackage{graphicx}      
\usepackage{natbib}        
\usepackage{todonotes}
\usepackage{times}
\usepackage{preamble}
\usepackage{tikz}
\usepackage{graphicx}
\usetikzlibrary{shapes,arrows}
\usepackage{pgfplots}
\pgfplotsset{compat=newest}
\usepackage{amsmath, amsfonts}
\usepackage[noend]{algpseudocode}
\usepackage{algorithm}
\acsetup{patch/maketitle = false}

\usepackage{dsfont}
\DeclareAcronym{aic}{
    short = AIC,
    long = Akaike Information Criterion 
}

\DeclareAcronym{aicc}{
    short = AICc,
    long = Akaike Information Criterion correction 
}

\DeclareAcronym{bfgs}{
    short = BFGS,
    long = Broyden-Fletcher-Goldfarb-Shanno algorithm
}

\DeclareAcronym{bic}{
    short = BIC,
    long = Bayesian Information Criterion 
}

\DeclareAcronym{cks}{
    short = CKS,
    long = Compositional Kernel Search 
}

\DeclareAcronym{dl}{
    short = DL,
    long = Deep Learning 
}

\DeclareAcronym{gp}{
    short = GP,
    long = Gaussian Process,
    long-plural-form = Gaussian Processes 
}

\DeclareAcronym{hmc}{
    short = HMC, 
    long = Hamiltonian Monte Carlo
}

\DeclareAcronym{kl}{
    short = KL, 
    long = Kullback-Leibler
}

\DeclareAcronym{ks}{
    short = KS,
    long = Kernel Search 
}

\DeclareAcronym{lfm}{
    short = LFM, 
    long = latent force model 
}

\DeclareAcronym{lodegp}{
    short = LODE-GP, 
    long = Linear Ordinary Differential Equation Gaussian Process,
    long-plural-form = Linear Ordinary Differential Equation Gaussian Processes
}

\DeclareAcronym{lti}{
    short = LTI, 
    long = Linear Time Invariant 
}

\DeclareAcronym{map}{
    short = MAP, 
    long = maximum a posteriori 
}

\DeclareAcronym{mc}{
    short = MC, 
    long = Monte Carlo 
}

\DeclareAcronym{mcmc}{
    short = MCMC, 
    long = Markov Chain Monte Carlo 
}

\DeclareAcronym{mll}{
    short = MLL, 
    long = marginal log likelihood 
}

\DeclareAcronym{ml}{
    short = ML, 
    long = Machine Learning
}

\DeclareAcronym{mpc}{
    short = MPC, 
    long = model predictive control 
}

\DeclareAcronym{lmpc}{
    short = LMPC, 
    long = linear model predictive control 
}

\DeclareAcronym{nmpc}{
    short = NMPC, 
    long = nonlinear model predictive control 
}

\DeclareAcronym{nn}{
    short = NN, 
    long = Neural Network 
}

\DeclareAcronym{nuts}{
    short = NUTS, 
    long = No U-Turn Sampler 
}

\DeclareAcronym{ode}{
    short = ODE, 
    long = ordinary differential equation 
}

\DeclareAcronym{rmse}{
    short = RMSE, 
    long = root mean squared error 
}

\DeclareAcronym{se}{
    short = SE, 
    long = squared exponential 
}

\DeclareAcronym{sgd}{
    short = SGD, 
    long = Stochastic Gradient Descent 
}

\DeclareAcronym{skc}{
    short = SKC, 
    long = Structured Kernel Composition 
}

\DeclareAcronym{snf}{
    short = SNF, 
    long = Smith Normal Form 
}

\DeclareAcronym{rkhs}{
    short = RKHS, 
    long = Reproducing Kernel Hilbert Space 
}

\begin{document}
\begin{frontmatter}

\title{Physics-informed Gaussian Processes for Model Predictive Control of Nonlinear Systems \thanksref{footnoteinfo}} %
\thanks[footnoteinfo]{Jörn Tebbe and Andreas Besginow are supported by the SAIL project which is funded by the Ministry of Culture and Science of the State of North Rhine-Westphalia under the grant no NW21-059C.}

\author[First]{Adrian Lepp} 
\author[Second]{Jörn Tebbe} 
\author[Second]{Andreas Besginow}

\address[First]{
Bielefeld University, 
CITEC
(e-mail: adrian.lepp@uni-bielefeld.de).}
\address[Second]{Institute Industrial IT - inIT, OWL University of Applied Sciences and Arts Lemgo, Germany (e-mail: joern.tebbe@th-owl.de)}

\begin{abstract}
Recently, a novel linear model predictive control algorithm based on a physics-informed Gaussian Process has been introduced, whose realizations strictly follow a system of underlying linear ordinary differential equations with constant coefficients. The control task is formulated as an inference problem by conditioning the Gaussian process prior on the setpoints and incorporating pointwise soft-constraints as further virtual setpoints. We apply this method to systems of nonlinear differential equations, obtaining a local approximation through the linearization around an equilibrium point. In the case of an asymptotically stable equilibrium point convergence is given through the Bayesian inference schema of the Gaussian Process. Results for this are demonstrated in a numerical example.
\end{abstract}


\end{frontmatter}

\section{Introduction}
Model predictive control (\acs{mpc}) is an advanced control method that is commonly applied in industrial applications for a wide range of
dynamic systems by reformulating the control task as an optimization problem \citep{rawlings2017model, tebbe2023holistic}.
It addresses the tracking problem, where the system states should follow a given reference signal \citep{Grüne_Pannek_2017}.
The basic principle of \ac{mpc} is to simulate the future behavior of the system with a \emph{predictive model} embedded in a \emph{control strategy}, which optimizes the control inputs with respect to an objective function and additional constraints.

While the optimization problem is a convex quadratic problem for linear systems, it is generally no longer convex for Nonlinear \ac{mpc} (NMPC) and finding the global optimum is not guaranteed.
To overcome this problem, several implementations for nonlinear systems exist that are based on their prior linearization, designing the controller for a linear surrogate model \citep{Zheng_2000, torrisi_projected_2016}. 
Linearization methods for dynamic systems can be global but more complex, like input-output linearization \citep{Kouvaritakis_Cannon_Rossiter_2000}, while others create multiple local models by performing successive linearizations \citep{Qin_Badgwell_2000}.

Classical predictive models are \emph{first-principle-based} but with the rise of machine learning, many \emph{data-driven} models have emerged \citep{draeger_model_1995, piche_nonlinear_2000, berberich_data-driven_2021}. Especially \acp{gp} are commonly applied in the modeling of dynamic systems due to their excellent handling of limited data and uncertainty quantification \citep{
hewing2018cautious, maiworm2021online}. 
If knowledge about the system dynamics exists, physics-informed \acp{gp} offer a combination of data-driven and first-principle-based methods in the form of assumptions about general system behavior \citep{pmlr-v5-alvarez09a, Ross_Smith_Alvarez_2021} or by directly incorporating differential equations 
\citep{besginow2022constraining, harkonen2023gaussian}.
While \ac{gp} based tracking \ac{mpc} schemes usually  incorporate the \ac{gp} only as a predictive model \citep{Kocijan_Murray-Smith_Rasmussen_Girard_2004, Umlauft_Beckers_Hirche_2018, Matschek_Himmel_Sundmacher_Findeisen_2020},  \cite{tebbe2024physicsinformedgaussianprocesseslinear} recently introduced a novel approach, that optimizes over the union of system dynamics and the control law in one \ac{gp} model.
To this end, they utilize the \ac{lodegp} \citep{besginow2022constraining}, a class of \acp{gp} that strictly satisfy given linear \ac{ode} systems.
This reduces the control task to a simple inference problem, also known as control as inference (CAI) which has been used in stochastic optimal control and reinforcement learning \citep{levine2018reinforcement}.
The optimal control problem is solved by incorporating the setpoints and constraints as training data for the \ac{lodegp} and obtaining the control law directly from its posterior predictive distribution.
Since the \ac{lodegp} does not distinguish inputs, state, and outputs, it is therefore a behavioral approach to control \citep{willems1997introduction}.

We extend this method to nonlinear systems by providing a linearization around an equilibrium point, which we also use as reference points for the control task. 
This equilibrium point is a state of the system that does not change over time. In many applications, we find the goals of steering the system towards this equilibrium point and stabilize the system there.
\cite{tebbe2024physicsinformedgaussianprocesseslinear} show that the kernelized structure of the \acp{lodegp} provides open-loop stability to the controlled system; this property also holds for nonlinear systems with an asymptotically stable equilibrium point as reference.
Including the reference as so-called equilibrium endpoint constraint yields finite-time convergence \citep{Grüne_Pannek_2017}.
Computation of linear surrogate models and the construction of the \ac{lodegp} can be done with  computer algebra \citep{oberst1990multidimensional,pommaret1999algebraic,zerz2000topics,chyzak2005effective,lange2013thomas,lange2020thomas}, allowing for automatic controller design.
We provide numerical results for a nonlinear two-tank system as an example.

\section{Problem formulation}
\label{sec:problem}
Consider the system of nonlinear \ac{ode}s
\begin{equation}\label{eq:nonlinear_system}
    \dot{\state}(t) = f(\state(t), \control(t))
\end{equation}
where  $\state(t) \in \mathbb{R}^{\dimState}$ contains the internal system states and $\control (t)\in \mathbb{R}^{\dimControl}$ is the control input.
In tracking control, it is the goal to find a control trajectory $\control(t)$ such that the system states $\state(t)$ follow a reference $\state_\text{ref}$ for a given initial state $\state(t_0) = \state_0$, while satisfying the state and control constraints 
\begin{align}\label{eq:box_constraints_x}
    \state_{\min} &\leq \state(t) \leq \state_{\max} \\
    \control_{\min} &\leq \control(t) \leq \control_{\max}\label{eq:box_constraints_u}
\end{align}
for $t \in [t_0, t_T]$. 
This task can be formulated as the minimization of a defined norm between the reference and states
\begin{equation}\label{eq:integral_formulation}
    \min_{\control(t)} \int_{t_0}^{t_T} \Vert \state_{\text{ref}} - \state(t) \Vert dt.
\end{equation}

Since control is often performed at discrete time steps, the tracking control task is reformulated as an approximation of \eqref{eq:integral_formulation}  to find the minimal error control solution using
\begin{subequations}
\begin{align} 
    \min_{\control(t)} &\sum_{i=0}^{T} \Vert \state_{\text{ref}} - \state(t_i) \Vert + \Vert \control(t) \Vert\label{eq:mpc:obj}\\
    \text{s.t. } \dot \state &= f(\state(t),\control(t)), \label{eq:mpc:con:ode}\\
    \state(t_0) &= \state_0, \label{eq:mpc:con:init}\\
    \state_{\min} &\leq \state(t) \leq \state_{\max} \quad \forall t \in [t_0,t_T], \label{eq:mpc:con:\state}\\
    \control_{\min} &\leq \control(t) \leq \control_{\max} \quad \forall t \in [t_0,t_T] \label{eq:mpc:con:u}.
\end{align}
\end{subequations}

This optimization problem is solved recursively for discrete timesteps with a moving horizon $t_T$, where the first element of the optimal control input $\control(t_0)$ is applied to the system for the next timestep.

While \ac{mpc} can be implemented for nonlinear systems, there exist several implementations, that linearize the system dynamics to design the controller for a linear surrogate model \citep{Zheng_2000}. The \ac{lodegp}-based \ac{mpc} algorithm presented by \cite{tebbe2024physicsinformedgaussianprocesseslinear} requires a system in the linear state-space form
\begin{align}\label{eq:linear_system}
    \dot \state(t) &= A \state(t) + B\control(t)
\end{align}
where $A \in \mathbb{R}^{\dimState \times \dimState}, B \in \mathbb{R}^{\dimState \times \dimControl}$. Thus, we have to find linear approximations of the nonlinear system \eqref{eq:mpc:con:ode}. We do this by linearizing the system around an equilibrium point $\state_e$, which also serves as the reference point $\state_{\text{ref}}$ to which we want to steer the system states.

\section{Preliminaries}
\label{sec:preliminaries}
\subsection{Linearization Around an Equilibrium Point}\label{sec:linearization}
In the following, we describe the linearization of nonlinear systems around equilibrium points by first introducing the concept of an equilibrium point according to \cite{Adamy_2022} and then showing how one can linearize a system around such a point. This is a common technique in control theory, as it allows us to find local linear approximations.
\begin{definition}
    Given the system $\dot \state(t) = f(\state(t),\control(t))$, an \emph{equilibrium point} ${\state}_e$ is a point in the state space that satisfies
    \begin{equation}\label{eq:equilibrium_point}
        \dot \state(t) = f({\state}_e, {\control}_e) = 0,
    \end{equation}
    where ${\control}_e$ is an arbitrary but constant control input.
\end{definition}
It is often desired to transfer the system's states to such an equilibrium point that remains constant over time and to hold it there.  
To this end, it is useful to examine the notion of asymptotic stability.
\begin{definition}An equilibrium point is \emph{asymptotically stable} if, for every $\epsilon$-neighborhood there exists a $\delta$-neighborhood such that every trajectory starting in the $\delta$-neighborhood stays in the $\epsilon$-neighborhood for all $t > 0$ and furthermore converges to the equilibrium $\state_e$ as
\begin{equation}
    \lim_{t \to \infty}||\state(t)-\state_e ||=0.
\end{equation}    
The $\delta$-neighborhood is then called the \emph{basin of attraction} of $\state_e$.    
\end{definition}
In the following, we will not discuss how to investigate asymptotic stability but assume that we can determine such equilibrium points by solving \eqref{eq:equilibrium_point} for ${\state}_e$ by choosing an appropriate control signal ${\control}_e$.
By defining the new delta coordinates $\Delta\state = \state - {\state}_e$ and $\Delta\control = \control - {\control}_e$ in the neighborhood of the equilibrium point, we rewrite the system equations in the form
\begin{equation}\label{eq:delta_coordinates}
    {\Delta\dot\state}(t) = \mathbf{f}(\Delta\state(t) + {\state}_e, \Delta\control(t) + {\control}_e).
\end{equation}
Now, we linearize the system using the Taylor expansion of first-order
\begin{equation}
    {\Delta\dot\state}(t) \approx 
    \mathbf{f}({\state}_e, {\control}_e) + 
    \frac{\partial \syseq}{\partial \state} \vert_{{\state}_e, {\control}_e}   \cdot \Delta\state(t) + 
    \frac{\partial \syseq}{\partial \control} \vert_{{\state}_e, {\control}_e} \cdot \Delta\control(t),
\end{equation}
and with $\mathbf{f}({\state}_e, {\control}_e) = \mathbf{0}$, the linearized state space model is given by
\begin{equation}\label{eq:linearStateSpace}
\begin{alignedat}{4}
    {\Delta\dot\state}(t) &\approx 
    \frac{\partial \syseq}{\partial \state} \vert_{{\state}_e, {\control}_e}   && \cdot \Delta\state(t)  && + 
    \frac{\partial \syseq}{\partial \control} \vert_{{\state}_e, {\control}_e} && \cdot \Delta\control(t) \\
&\doteq {A}_e  && \cdot \Delta\state(t) &&+ {B}_e && \cdot \Delta\control(t)
\end{alignedat}
\end{equation}
with the constant Jacobian matrices ${A}_e, {B}_e$.
Equation~\eqref{eq:linearStateSpace} approximates the nonlinear form for small deviations around the equilibrium point. 

\subsection{Gaussian Processes}

A \acf{gp} \citep{rasmussen2006gaussian} 
\begin{equation}
    g(t) \sim \mathcal{GP}(\mu(t), k(t, t'))    
\end{equation}
is a stochastic process with the property that all random variables $g(t_1), \ldots, g(t_n)$ follow a jointly Gaussian distribution.
It is fully characterized by its mean function
\begin{equation}
    \mu(t):=\mathbb{E}[g(t)]
\end{equation}
and
covariance function
\begin{equation}\label{eq:covariance}
    k(t, t') := \mathbb{E}[(g(t) - \mu(t))(g(t') - \mu(t'))^{\top}]
    .
\end{equation}
By conditioning a \ac{gp} on a noisy dataset
$$\mathcal{D} = \{(t_1, z_1), \ldots ,(t_n, z_n)\}$$ 
with $z \sim g(t) + \mathcal{N}(0, \sigma_n^2)$
we can obtain the posterior \ac{gp}
\begin{equation}\label{eq:gaussian_process_posterior_distribution}
    \begin{aligned}
        \mu^* &= \mu(t^*) + K_*^T(K + \sigma_n^2 I)^{-1}z \\
        k^* &= K_{**} - K_*^T(K + \sigma_n^2 I)^{-1} {K_*}
    \end{aligned}
\end{equation}
with covariance matrices $K = (k(t_i, t_j))_{i,j} \in \mathbb{R}^{n \times n}$, $K_* = (k(t_i, t^*_j))_{i, j} \in \mathbb{R}^{n \times m}$ and $K_{**} = (k(t^*_i, t^*_j))_{i, j} \in \mathbb{R}^{m \times m}$ for predictive positions $t^* \in \mathbb{R}^m$ with noise variance $\sigma_n^2$.

While the mean function is often chosen as $\mu(t) = 0$ \citep{rasmussen2006gaussian}, a popular choice for the covariance function is the \acf{se} kernel
\begin{equation}\label{eq:SE_kernel}
    k_{\text{SE}}(t, t') = \sigma_f^2\exp\left(-\frac{(t-t')^2}{2\ell^2}\right)
    ,
\end{equation}
assuming smooth and infinitely differentiable functions.
Together with the noise variance $\sigma_n^2$, the signal variance $\sigma_f^2$ and lengthscale $\ell^2$ define a set of hyperparameters $\theta$, which are trained by maximizing the \ac{gp} \ac{mll}
\begin{equation} \label{eq:gp:MLL}
    \log p(z|t) = -\frac{1}{2}z^T\left(K_z \right)^{-1}z - \frac{1}{2}\log \left( \det \left(K_z \right) \right)
\end{equation}
where $K_z=K+\sigma_n^2I$,
$I$ is the identity and constant terms are omitted. We obtain a quadratic type error term combined with a regularization term based on the determinant of the regularized kernel matrix.

In this work, we will also exploit the two following properties of \acp{gp}:
First, the observation noise $\sigma_n^2$ can be input-dependent (heteroscedastic), i.e. $\sigma_n^2(t) \in \mathbb{R}^{n}$, allowing us to set individual noise levels for different datapoints.
Second, it is possible to manipulate existing \acp{gp} by applying a linear operator $\mathcal{L}$ on a \ac{gp} $g(t) \sim \mathcal{GP} (\mu(t), k(t,t^\prime)$, leading to another \ac{gp}
\begin{equation}\label{eq:gp:linear_operator}
    \mathcal{L}g(t) \sim \mathcal{GP} (\mathcal{L} \mu(t), \mathcal{L}k(t,t^\prime)\mathcal{L}^{\prime^\top}),   
\end{equation}
where $\mathcal{L}^{\prime}$ is the application of $\mathcal{L}$ on $t^\prime$
\citep{jidling2017linearly,langehegermann2018algorithmic}.

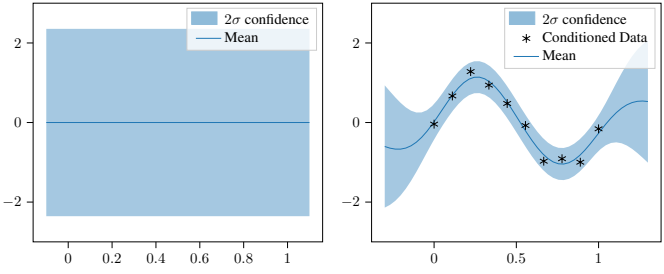
\begin{figure}\label{fig:GP_example}
    \centering
    \scalebox{0.97}{
\begin{tikzpicture}[scale=0.573]

\definecolor{darkgray176}{RGB}{176,176,176}
\definecolor{lightgray204}{RGB}{204,204,204}
\definecolor{steelblue31119180}{RGB}{31,119,180}

\definecolor{tab10_blue}{RGB}{31,119,180}

\begin{axis}[
legend cell align={left},
legend style={fill opacity=0.8, draw opacity=1, text opacity=1, draw=lightgray204},
tick align=outside,
tick pos=left,
x grid style={darkgray176},
xmin=-0.160000002756715, xmax=1.16000002510846,
xtick style={color=black},
y grid style={darkgray176},
ymin=-3, ymax=3,
ytick style={color=black}
]
\path [fill=steelblue31119180, fill opacity=0.4]
(axis cs:-0.100000001490116,2.35490489006042)
--(axis cs:-0.100000001490116,-2.35490489006042)
--(axis cs:-0.0142857134342194,-2.35490489006042)
--(axis cs:0.0714285746216774,-2.35490489006042)
--(axis cs:0.157142877578735,-2.35490489006042)
--(axis cs:0.242857158184052,-2.35490489006042)
--(axis cs:0.328571438789368,-2.35490489006042)
--(axis cs:0.414285749197006,-2.35490489006042)
--(axis cs:0.5,-2.35490489006042)
--(axis cs:0.585714280605316,-2.35490489006042)
--(axis cs:0.671428561210632,-2.35490489006042)
--(axis cs:0.757142901420593,-2.35490489006042)
--(axis cs:0.842857122421265,-2.35490489006042)
--(axis cs:0.928571462631226,-2.35490489006042)
--(axis cs:1.0142856836319,-2.35490489006042)
--(axis cs:1.10000002384186,-2.35490489006042)
--(axis cs:1.10000002384186,2.35490489006042)
--(axis cs:1.10000002384186,2.35490489006042)
--(axis cs:1.0142856836319,2.35490489006042)
--(axis cs:0.928571462631226,2.35490489006042)
--(axis cs:0.842857122421265,2.35490489006042)
--(axis cs:0.757142901420593,2.35490489006042)
--(axis cs:0.671428561210632,2.35490489006042)
--(axis cs:0.585714280605316,2.35490489006042)
--(axis cs:0.5,2.35490489006042)
--(axis cs:0.414285749197006,2.35490489006042)
--(axis cs:0.328571438789368,2.35490489006042)
--(axis cs:0.242857158184052,2.35490489006042)
--(axis cs:0.157142877578735,2.35490489006042)
--(axis cs:0.0714285746216774,2.35490489006042)
--(axis cs:-0.0142857134342194,2.35490489006042)
--(axis cs:-0.100000001490116,2.35490489006042)
--cycle;

\addlegendentry{$2\sigma$ confidence}
\addlegendimage{
      confidence box={circle,fill,inner sep=1pt}{diamond,fill,inner sep=1pt},
      steelblue31119180}

\addplot [semithick, steelblue31119180]
table {%
-0.100000001490116 0
-0.0142857134342194 0
0.0714285746216774 0
0.157142877578735 0
0.242857158184052 0
0.328571438789368 0
0.414285749197006 0
0.5 0
0.585714280605316 0
0.671428561210632 0
0.757142901420593 0
0.842857122421265 0
0.928571462631226 0
1.0142856836319 0
1.10000002384186 0
};
\addlegendentry{Mean}
\end{axis}

\end{tikzpicture}}
    \scalebox{0.97}{
\begin{tikzpicture}[scale=0.573]

\definecolor{darkgray176}{RGB}{176,176,176}
\definecolor{lightgray204}{RGB}{204,204,204}
\definecolor{steelblue31119180}{RGB}{31,119,180}

\begin{axis}[
legend cell align={left},
legend style={fill opacity=0.8, draw opacity=1, text opacity=1, draw=lightgray204},
tick align=outside,
tick pos=left,
x grid style={darkgray176},
xmin=-0.38000001013279, xmax=1.37999995052814,
xtick style={color=black},
y grid style={darkgray176},
ymin=-3, ymax=3,
ytick style={color=black}
]
\path [fill=steelblue31119180, fill opacity=0.4]
(axis cs:-0.300000011920929,0.934599936008453)
--(axis cs:-0.300000011920929,-2.13949489593506)
--(axis cs:-0.26800000667572,-2.07047033309937)
--(axis cs:-0.236000016331673,-1.96931564807892)
--(axis cs:-0.204000025987625,-1.83424127101898)
--(axis cs:-0.172000020742416,-1.66563546657562)
--(axis cs:-0.140000015497208,-1.46646726131439)
--(axis cs:-0.10800002515316,-1.24252831935883)
--(axis cs:-0.0760000348091125,-1.00236678123474)
--(axis cs:-0.0440000295639038,-0.756528854370117)
--(axis cs:-0.0120000317692757,-0.515527546405792)
--(axis cs:0.0199999660253525,-0.286729007959366)
--(axis cs:0.0519999638199806,-0.0727817714214325)
--(axis cs:0.0839999616146088,0.125859975814819)
--(axis cs:0.115999966859818,0.307051777839661)
--(axis cs:0.147999957203865,0.46556892991066)
--(axis cs:0.179999947547913,0.594264984130859)
--(axis cs:0.211999952793121,0.685999631881714)
--(axis cs:0.24399995803833,0.73531174659729)
--(axis cs:0.275999963283539,0.739395380020142)
--(axis cs:0.307999938726425,0.698349714279175)
--(axis cs:0.339999943971634,0.614844381809235)
--(axis cs:0.371999949216843,0.493481397628784)
--(axis cs:0.403999924659729,0.340139895677567)
--(axis cs:0.435999929904938,0.161466330289841)
--(axis cs:0.467999935150146,-0.0354481935501099)
--(axis cs:0.499999940395355,-0.24323758482933)
--(axis cs:0.531999945640564,-0.454380929470062)
--(axis cs:0.563999950885773,-0.661306738853455)
--(axis cs:0.595999956130981,-0.856574296951294)
--(axis cs:0.627999901771545,-1.03311252593994)
--(axis cs:0.659999907016754,-1.1845269203186)
--(axis cs:0.691999912261963,-1.30537581443787)
--(axis cs:0.723999977111816,-1.3913506269455)
--(axis cs:0.755999982357025,-1.43936347961426)
--(axis cs:0.787999987602234,-1.44754993915558)
--(axis cs:0.819999992847443,-1.41537308692932)
--(axis cs:0.851999998092651,-1.3440535068512)
--(axis cs:0.883999943733215,-1.23734652996063)
--(axis cs:0.915999948978424,-1.10260486602783)
--(axis cs:0.947999954223633,-0.951777815818787)
--(axis cs:0.979999959468842,-0.801335334777832)
--(axis cs:1.01199996471405,-0.669523298740387)
--(axis cs:1.04399991035461,-0.571118831634521)
--(axis cs:1.07599997520447,-0.513702988624573)
--(axis cs:1.10799992084503,-0.498257458209991)
--(axis cs:1.13999998569489,-0.521744549274445)
--(axis cs:1.17199993133545,-0.57893979549408)
--(axis cs:1.2039999961853,-0.663384795188904)
--(axis cs:1.23599994182587,-0.768006265163422)
--(axis cs:1.26800000667572,-0.885636508464813)
--(axis cs:1.29999995231628,-1.00948870182037)
--(axis cs:1.29999995231628,2.06460523605347)
--(axis cs:1.29999995231628,2.06460523605347)
--(axis cs:1.26800000667572,1.96687960624695)
--(axis cs:1.23599994182587,1.83457255363464)
--(axis cs:1.2039999961853,1.66752183437347)
--(axis cs:1.17199993133545,1.46814382076263)
--(axis cs:1.13999998569489,1.24169707298279)
--(axis cs:1.10799992084503,0.996329128742218)
--(axis cs:1.07599997520447,0.742780327796936)
--(axis cs:1.04399991035461,0.493383884429932)
--(axis cs:1.01199996471405,0.259850561618805)
--(axis cs:0.979999959468842,0.049962043762207)
--(axis cs:0.947999954223633,-0.134039491415024)
--(axis cs:0.915999948978424,-0.292916357517242)
--(axis cs:0.883999943733215,-0.426429897546768)
--(axis cs:0.851999998092651,-0.531766712665558)
--(axis cs:0.819999992847443,-0.604436576366425)
--(axis cs:0.787999987602234,-0.639964938163757)
--(axis cs:0.755999982357025,-0.635309278964996)
--(axis cs:0.723999977111816,-0.589667439460754)
--(axis cs:0.691999912261963,-0.504617929458618)
--(axis cs:0.659999907016754,-0.383751332759857)
--(axis cs:0.627999901771545,-0.232107102870941)
--(axis cs:0.595999956130981,-0.0556027591228485)
--(axis cs:0.563999950885773,0.13930931687355)
--(axis cs:0.531999945640564,0.345812737941742)
--(axis cs:0.499999940395355,0.556771397590637)
--(axis cs:0.467999935150146,0.764747262001038)
--(axis cs:0.435999929904938,0.962083578109741)
--(axis cs:0.403999924659729,1.1411120891571)
--(axis cs:0.371999949216843,1.29448628425598)
--(axis cs:0.339999943971634,1.41562008857727)
--(axis cs:0.307999938726425,1.49910879135132)
--(axis cs:0.275999963283539,1.54107928276062)
--(axis cs:0.24399995803833,1.53936886787415)
--(axis cs:0.211999952793121,1.49358296394348)
--(axis cs:0.179999947547913,1.40519857406616)
--(axis cs:0.147999957203865,1.27785742282867)
--(axis cs:0.115999966859818,1.1179713010788)
--(axis cs:0.0839999616146088,0.935551404953003)
--(axis cs:0.0519999638199806,0.744957208633423)
--(axis cs:0.0199999660253525,0.5645672082901)
--(axis cs:-0.0120000317692757,0.413845837116241)
--(axis cs:-0.0440000295639038,0.307981014251709)
--(axis cs:-0.0760000348091125,0.254113435745239)
--(axis cs:-0.10800002515316,0.252056956291199)
--(axis cs:-0.140000015497208,0.296977400779724)
--(axis cs:-0.172000020742416,0.381446242332458)
--(axis cs:-0.204000025987625,0.496665835380554)
--(axis cs:-0.236000016331673,0.633264183998108)
--(axis cs:-0.26800000667572,0.782045543193817)
--(axis cs:-0.300000011920929,0.934599936008453)
--cycle;
\addlegendentry{$2\sigma$ confidence}
\addlegendimage{
      confidence box={circle,fill,inner sep=1pt}{diamond,fill,inner sep=1pt},
      steelblue31119180}

\addplot [semithick, black, mark=asterisk, mark size=3, mark options={solid}, only marks]
table {%
0 -0.0409419648349285
0.111111111938953 0.667891800403595
0.222222223877907 1.28061735630035
0.333333343267441 0.940110445022583
0.444444447755814 0.478916764259338
0.555555582046509 -0.0743785798549652
0.666666626930237 -0.972533524036407
0.777777791023254 -0.910220146179199
0.888888895511627 -0.997552394866943
1 -0.158487573266029
};
\addlegendentry{Conditioned Data}
\addplot [semithick, steelblue31119180]
table {%
-0.300000011920929 -0.60244745016098
-0.26800000667572 -0.644212424755096
-0.236000016331673 -0.668025732040405
-0.204000025987625 -0.668787717819214
-0.172000020742416 -0.642094612121582
-0.140000015497208 -0.584744930267334
-0.10800002515316 -0.495235681533813
-0.0760000348091125 -0.374126672744751
-0.0440000295639038 -0.224273920059204
-0.0120000317692757 -0.0508408546447754
0.0199999660253525 0.138919115066528
0.0519999638199806 0.336087703704834
0.0839999616146088 0.530705690383911
0.115999966859818 0.712511539459229
0.147999957203865 0.871713161468506
0.179999947547913 0.999731779098511
0.211999952793121 1.0897912979126
0.24399995803833 1.13734030723572
0.275999963283539 1.14023733139038
0.307999938726425 1.09872925281525
0.339999943971634 1.01523220539093
0.371999949216843 0.893983840942383
0.403999924659729 0.740625977516174
0.435999929904938 0.561774969100952
0.467999935150146 0.364649534225464
0.499999940395355 0.156766891479492
0.531999945640564 -0.0542840957641602
0.563999950885773 -0.260998725891113
0.595999956130981 -0.456088542938232
0.627999901771545 -0.632609844207764
0.659999907016754 -0.784139156341553
0.691999912261963 -0.904996871948242
0.723999977111816 -0.990509033203125
0.755999982357025 -1.0373363494873
0.787999987602234 -1.04375743865967
0.819999992847443 -1.0099048614502
0.851999998092651 -0.937910079956055
0.883999943733215 -0.831888198852539
0.915999948978424 -0.697760581970215
0.947999954223633 -0.542908668518066
0.979999959468842 -0.375686645507812
1.01199996471405 -0.204836368560791
1.04399991035461 -0.0388674736022949
1.07599997520447 0.114538669586182
1.10799992084503 0.249035835266113
1.13999998569489 0.359976291656494
1.17199993133545 0.444602012634277
1.2039999961853 0.502068519592285
1.23599994182587 0.533283174037933
1.26800000667572 0.54062157869339
1.29999995231628 0.527558207511902
};
\addlegendentry{Mean}
\end{axis}

\end{tikzpicture}}
    \caption{(Left) A \ac{gp} prior with zero mean and \ac{se} covariance function. (Right) The same \ac{gp}, but conditioned on datapoints (black asterisk).
    The blue line is its mean and the blue area is two times its standard deviation ($2\sigma$).}
\end{figure}

\subsection{Linear Ordinary Differential Equation GPs}\label{sec:LODE_GP}
We review the construction of the \ac{lodegp} --- a \ac{gp} that strictly satisfies an underlying system of linear homogeneous ordinary differential equations --- as introduced in \citep{besginow2022constraining} by starting with a linearized system in general state space form
\begin{equation}
    {\Delta\dot\state}(t) = A_e \cdot  \Delta\state(t) + B_e \cdot \Delta\control(t)
\end{equation}
in delta coordinates. 
First, the system representation must be changed to a set of homogeneous differential equations by subtracting ${\Delta\dot\state}$ and stacking the state $\Delta\state$ and the input $\Delta\control$ in one variable $\Delta\stateZ$ with
\begin{equation}\label{eq:lodegp_homogeneous}
    0 = H \cdot \Delta\stateZ(t) = \left[
    \begin{array}{c|c}
         A_e - \bm{I} \cdot \partial_t & B_e
    \end{array}
    \right] 
    \cdot 
    \begin{bmatrix}
        {\Delta\state}(t) \\ {\Delta\control}(t)     
    \end{bmatrix}
\end{equation}
where  
$H = 
\left[ \begin{array}{c|c}
    A_e - \bm{I} \cdot  \partial_t & B_e
\end{array} \right]$
is a $\dimState \times \dimStateZ$ operator matrix with $\dimStateZ=\dimState + \dimControl$, $\bm{I}$ the identity matrix of size $\dimState \times \dimState$ and $\partial_t$ the differential operator.
The matrix $H$
can be decoupled by calculating the Smith Normal Form \citep{smith1862systems, newman1997smith}
\begin{equation}
    D = W \cdot H \cdot V,
\end{equation}
with diagonal matrix $D$ and invertible square matrices $W$ and $V$.
All matrices are operator matrices and thus belong to the polynomial ring $\mathbb{W}[\partial_t]$.
Left multiplication with $W$ and  neutral multiplication
with $V \cdot V^{-1}$ of Equation~\eqref{eq:lodegp_homogeneous} yields
\begin{align}\label{eq:lodegp40}
    W \cdot H \cdot V \cdot V^{-1} \Delta\stateZ(t) &= 0 \nonumber \\
    D \cdot V^{-1} \Delta\stateZ(t) &= 0
    .
\end{align}
Introducing the latent state vector 
\begin{equation}\label{eq:lodegp_latent}
\latent(t)=V^{-1} \Delta\stateZ(t)    
\end{equation}
allows, to rewrite the system with
\begin{equation}\label{eq:lodegp50}
    D \cdot \latent(t) = 0
\end{equation}
to obtain a decoupled system of linear ordinary differential equations with diagonal matrix $D$. This decoupling introduces two useful possibilities. First the independence of the single dimensions of $\latent$ allows constructing a $\dimStateZ$-dimensional latent \ac{gp}
for  $\latent(t)$
\begin{equation}\label{eq:latentgp}
    h(t)
    \sim \mathcal{GP}(0, k(t, t'))
    ,
\end{equation}
where $k(t, t')$ is a multidimensional covariance function with dimensionality $\dimStateZ$.
Furthermore, one can easily  determine independent solutions for $\latent(t)$. Together with 
the definition of the covariance function in Equation~\eqref{eq:covariance} and
the mean function $\mu(t)=0$
it is possible to construct a covariance function for the latent \ac{gp} $h(t)$ containing the solutions of the ODEs.
The entries of $D$ are either given with zero, one or a polynomial and \cite{besginow2022constraining} provide a set of rules to construct the covariance function without a need to actually solve the equations. In the case of zeros in the diagonal entries of $D$, this indicates a degree of freedom in the system, usually introduced by a control input. \cite{besginow2022constraining} propose to use an \ac{se} kernel for these entries, which allows adapting the degrees of freedom to given  data.

Applying the inverse transformation of Equation~\eqref{eq:lodegp_latent} on $h(t)$ following Equation~\eqref{eq:gp:linear_operator} yields
the \ac{lodegp} 
\begin{equation}
    \Delta g(t) =
    Vh(t)
    \sim
    \mathcal{GP}(0, V k(t,t') V'^{\top})
\end{equation}
over $[\Delta \state(t), \Delta \control(t)]^\top$.
Finally, using the transformation $\state= \Delta \state + \state_e$, $\control=\Delta \control + \control_e$ 
leads to another \ac{lodegp} 
\begin{equation}\label{eq:finalGP}
    g(t) = \Delta g(t) + 
    \begin{bmatrix}{\state}_e \\ {\control}_e\end{bmatrix}
    \sim
    \mathcal{GP}\left(
        \begin{bmatrix}{\state}_e \\ {\control}_e\end{bmatrix}, 
        V k(t,t') V'^{\top}
    \right),
\end{equation}
which outputs continuous solutions for $\state(t)$ and $\control(t)$ that approximately fulfill the underlying nonlinear differential equations in the neighborhood of $\state_e$.
Training and conditioning the \ac{lodegp} on datapoints furthermore adapts the included degrees of freedom to fit the data. Since these are encoded in the control input, this makes it possible to add datapoints as setpoints to the system and find control trajectories that lead to the desired behavior. This property of the \ac{gp} will be exploited in the next chapter to construct a \ac{mpc} algorithm.

\section{LODE-GP Model Predictive Control}\label{sec:method}
In this section we formulate the \ac{mpc} problem as a \ac{lodegp} inference problem as introduced by \cite{tebbe2024physicsinformedgaussianprocesseslinear}.
We want to consider the transition between two setpoints as a special case of reference as done by \cite{Matschek_Himmel_Sundmacher_Findeisen_2020}. In the transient phase 
the system should change from the initial state $\state_0$ to the constant reference $\state_\text{ref}$  in finite time $t_{\text{ref}}$ and stay constant in the asymptotic phase afterward. 

At this point we want to emphasize the difference from other \ac{gp}-based \ac{mpc} approaches as \citep{Kocijan_Murray-Smith_Rasmussen_Girard_2004, Matschek_Himmel_Sundmacher_Findeisen_2020}. These approaches use the \ac{gp} to directly model the system dynamics with $\dot{\state}(t)=f(\state(t), \control(t)) \sim \mathcal{GP}(0, k)$ and the system dynamics are incorporated completely via datapoints. The \ac{mpc} algorithm optimizes over $\control(t)$ and predicts the states with the \ac{gp} to find an optimal trajectory that fulfills the constraints. 

In contrast, the \ac{lodegp} already incorporates the system dynamics in the kernel function and is used as a predictor for the states and the control input together. Constraints have to be expressed in the form of datapoints
and the optimization over $\control(t)$ is obsolete since the \ac{lodegp} posterior adapts  the degrees of freedom in the form of the control input to satisfy the constraints.

In the following, we will first show how the constraints \eqref{eq:mpc:con:ode} and \eqref{eq:mpc:con:init} can be satisfied, before we discuss how open-loop stability can be guaranteed and how convergence to the desired setpoint in finite time can be achieved.

\subsection{Implementing Hard and Soft Constraints} \label{sec:method:basic}

We assume that we can linearize the system in Equation~\eqref{eq:nonlinear_system} around an asymptotically stable equilibrium point $\state_e$ and all state trajectories starting in $\state_0$ stay in its basin of attraction. For this approximation we can then construct a \ac{lodegp} over $\stateZ(t)=[\state(t), \control(t)]^\top$ 
from Equation~\eqref{eq:finalGP} 
which is a valid model of the system dynamics in the considered region and hence satisfies \eqref{eq:mpc:con:ode}.

The other constraints are incorporated pointwise in the dataset $\mathcal{D}$ into the \ac{lodegp} by using  heteroscedastic noise  $\sigma_n^2 \in \mathbb{R}_{\geq 0}^{\dimStateZ}$ allowing for different noise levels on each state and control dimension. 
We respect the initial point constraint~\eqref{eq:mpc:con:init} as a hard constraint, by conditioning the \ac{lodegp} on the current state $\state_0$ and control input $\control_0$ in the dataset $\mathcal{D}_\text{init} = \{(t_0, \stateZ_0)\}$  using zero noise variance $\sigma_n^2(t_i) = 0 \in \mathbb{R}^{\dimStateZ}$ at every timestep.
Due to numerical issues, we have to set a numerical jitter of $10^{-8}$ as the noise variance. 
This forces the \ac{lodegp} posterior mean to satisfy $\mu^*(t_0) = \stateZ_0$ up to numerical precision.

The state and control constraints \eqref{eq:mpc:con:x}-\eqref{eq:mpc:con:u} are encoded as pointwise soft constraints in the dataset 
$$\mathcal{D}_\text{con} = \{(t_{1}, {\stateZ}_\text{con}), \dots,(t_{m_c}, {\stateZ}_\text{con})\}$$ 
with
\begin{equation}\label{eq:soft_constraint}
    {\stateZ}_\text{con} = \frac{({\stateZ}_\text{max} + {\stateZ}_\text{min})}{2}
\end{equation}
and the constraint noise variance
\begin{equation}\label{eq:soft_constraint_var}
    \sigma_n^2 = \frac{({\stateZ}_\text{max} - {\stateZ}_\text{min})}{4}.
\end{equation}

Technically, the incorporation of $\mathcal{D}= {D}_\text{init} \cup  \mathcal{D}_\text{con}$ only imposes soft constraints in the likelihood in Equation~\eqref{eq:gp:MLL} and therefore in the posterior mean $\mu^*(t)$ of the \ac{lodegp}, but setting small noise variances $\sigma_n^2$ enforces $\mu^*(t)$ to match these datapoints with high probability.

\subsection{Open-Loop Stability and Convergence in Finite Time}
With the choice of $\state_\text{ref}=\state_e$, the equilibrium $\state_e$ being asymptotically stable and with the initial state $\state_0$ starting in its basin of attraction, simply applying the control input $\control_e$ to the system would drive its states in infinite time to the equilibrium, hence to the reference.
This property also holds for the \ac{gp} over $\stateZ(t)$ with prior $\mu(t)=\stateZ_\text{ref}=[\state_e, \control_e]^\top$ and thus for our controller, since the posterior mean $\mu^*(t)$ of the underlying \ac{lodegp} converges to its prior for $t \to \pm \infty$, as proven by \cite{tebbe2024physicsinformedgaussianprocesseslinear}.

To achieve convergence to $\state_\text{ref}$ in finite time $t_{\text{ref}}$, we can incorporate the reference as an endpoint constraint \citep{Grüne_Pannek_2017} in the dataset $\mathcal{D}$ with $\mathcal{D}_\text{ref} = \{(t_{\text{ref}}, \stateZ_\text{ref})\}$ .
The \ac{gp} will then be used to create a control trajectory which translates the system from the initial setpoint $\state_0$ to the endpoint $\state_\text{ref}$ in $t_{\text{ref}}$. For $t > t_{\text{ref}}$ it is then sufficient to set the control input to $\control_e$ to keep the system constant at the reference. 
Note that it is not guaranteed to reach any endpoint $\state_\text{ref}$ in time $t_{\text{ref}}$ and choosing both values unreasonably may lead to violation of the constraints. Thus, one needs to respect the system dynamics and constraints, which is a trade-off between fast convergence and low overshoots.

\section{Evaluation}\label{sec:evaluation}

\subsection{System Description}
We consider a nonlinear system consisting of two water tanks as presented in Figure \ref{fig:watertank}. The system's behavior follows the nonlinear state space equations
\begin{align}\label{5_dreitankEq}
\dot{x}_1(t)
&=
\frac{1}{A} \left( u_1(t) - Q  \right )
\nonumber
, \\
\dot{x}_2(t)
&=
\frac{1}{A}\left ( Q
-c_{2\text{R}}\sqrt{2g x_2} \right)
, \\ 
Q &= c_{12}\cdot \text{sign}(x_1(t)-x_2(t))\sqrt{2g\mid x_1(t)-x_2(t)\mid}, 
\nonumber
\end{align}
where $x_1(t), x_2(t)$ represent the water levels in the tanks and $u_1(t)$ is the control input. A description of all parameters is given in Table ~\ref{tab:parameters}.
\definecolor{tab10_blue}{RGB}{31,119,180}

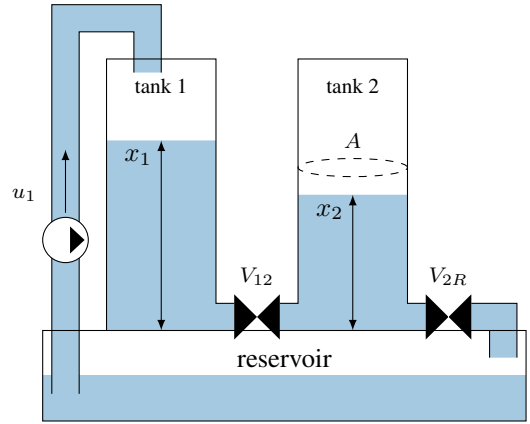
\begin{figure}[h!t]
    \centering
        
    \begin{tikzpicture}
        \def \factor{1.2}
        \def \noTank{2}

        \def \hT{\factor*3}
        \def \hR{\factor*1}
        
        \def \bT{\factor*1.2} 
        \def \dT{\factor*0.9}   
        \def \d{-\factor*0.2}    
        \def \bR{5.3*\factor}     
        
        \def \rohr{\factor*0.3} 
        \def \v{\factor*0.5}
        \def \a{\factor*0.4} 
        \def \b{\factor*0.1}
        \def \rV{\a *0.5} 

        \fill[tab10_blue,  opacity=0.4] ({\d+\noTank*\dT+\noTank*\bT+\rV+\v},{\hR+\rohr}) -- ({\d+\noTank*\dT+\noTank*\bT+\rV+\v+\v},{\hR+\rohr})--
        ({\d+\noTank*\dT+\noTank*\bT+\rV+\v+\v},{\hR-\rohr}) --
        ({\d+\noTank*\dT+\noTank*\bT+\rV+\v+\v-\rohr},{\hR-\rohr}) --
        ({\d+\noTank*\dT+\noTank*\bT+\rV+\v+\v-\rohr},{\hR}) --
        ({\d+\noTank*\dT+\noTank*\bT+\v+\v-\rohr},{\hR}) --
        ({\d+\noTank*\dT+\noTank*\bT+\v+\v-\rohr},{\hR+\rohr});
        \fill[tab10_blue!40!white]    ({\dT +\d},\hR) -- 
                ({\dT +\d},{\hR+\hT*0.7}) --
                ({\dT +\d + \bT},{\hR+\hT*0.7}) -- 
                ({\dT +\d + \bT},{\hR+\rohr}) --
                ({\dT +\d + \bT +\rV},{\hR+\rohr}) --
                ({\dT +\d + \bT +\rV},{\hR}) --
                ({\dT +\d},\hR)
                ;
        \fill[tab10_blue!40!white]    ({\d+\dT+\bT+\dT-\rV},{\hR+\rohr}) --
                ({\d+\dT+\bT+\dT},{\hR+\rohr}) --
                ({\d+\dT+\bT+\dT},{\hR+ \hT*0.5}) --
                ({\d+\dT+\bT+\dT +\bT},{\hR+\hT*0.5}) --
                ({\d+\dT+\bT+\dT +\bT},{\hR+\rohr}) --
                ({\d+\dT+\bT+\dT +\bT+\rV},{\hR+\rohr}) --
                ({\d+\dT+\bT+\dT +\bT+\rV},{\hR}) --
                ({\d+\dT+\bT+\dT -\rV},{\hR}) --
                ({\d+\dT+\bT+\dT-\rV},{\hR+\rohr})
                ;
        \fill[tab10_blue!40!white] (0,0) rectangle ++(\bR ,{\hR*0.5});
        \fill[tab10_blue!40!white]     ({\b},{\rohr}) --
                    ({\b},{\rohr*2+\hR+\hT}) --
                    ({\b+\dT+\rohr},{\rohr*2+\hR+\hT}) --
                    ({\b+\dT+\rohr},{\hR+\hT-\rohr*0.5})--
                    ({\b+\dT},{\hR+\hT-\rohr*0.5}) --
                    ({\b+\dT},{\rohr+\hR+\hT}) --
                    ({\b+\rohr},{\rohr+\hR+\hT}) --
                    ({\b+\rohr},{\rohr}) --
                    ({\b},{\rohr});
        
        \draw[draw=black] (0,0) rectangle ++(\bR,\hR);
        \draw   ({\dT +\d},\hR) -- 
                ({\dT +\d},{\hR+\hT}) --
                ({\dT +\d + \bT},{\hR+\hT}) -- 
                ({\dT +\d + \bT},{\hR+\rohr}) --
                ({\dT +\d + \bT +\rV},{\hR+\rohr})
                ;
        \draw   ({\d+\dT+\bT+\dT-\rV},{\hR+\rohr}) --
                ({\d+\dT+\bT+\dT},{\hR+\rohr}) --
                ({\d+\dT+\bT+\dT},{\hR+ \hT}) --
                ({\d+\dT+\bT+\dT +\bT},{\hR+\hT}) --
                ({\d+\dT+\bT+\dT +\bT},{\hR+\rohr}) --
                ({\d+\dT+\bT+\dT +\bT+\rV},{\hR+\rohr})
                ;
        
        \fill[black]    ({\d+\dT+\bT+\rV},{\hR-\b}) --                                        ({\d+\dT+\bT+\rV},{\hR+\rohr+\b}) -- 
                        ({\d+\dT+\bT+\rV+\v},{\hR-\b}) --
                        ({\d+\dT+\bT+\rV+\v},{\hR+\rohr+\b}) --
                        ({\d+\dT+\bT+\rV},{\hR-\b})
                        ;
        \node at ({\d+\dT+\bT+\rV+\v*0.5},{\hR+\v+\b}) {\small $V_{12}$};
        \fill[black]    ({\d+\dT+\bT+\rV+\bT+\dT},{\hR-\b}) --                                ({\d+\dT+\bT+\rV+\bT+\dT},{\hR+\rohr+\b}) -- 
                        ({\d+\dT+\bT+\rV+\v+\bT+\dT},{\hR-\b}) --
                        ({\d+\dT+\bT+\rV+\v+\bT+\dT},{\hR+\rohr+\b}) --
                        ({\d+\dT+\bT+\rV+\bT+\dT},{\hR-\b});
        \node at ({\d+2*\dT+2*\bT+\rV+\v*0.5},{\hR+\v+\b}) 
        {\small $V_{2R}$};
        
        \draw[] ({\d+\noTank*\dT+\noTank*\bT+\rV+\v},{\hR+\rohr}) -- ({\d+\noTank*\dT+\noTank*\bT+\rV+\v+\v},{\hR+\rohr})--
        ({\d+\noTank*\dT+\noTank*\bT+\rV+\v+\v},{\hR-\rohr});
        \draw[] ({\d+\noTank*\dT+\noTank*\bT+\rV+\v+\v-\rohr},{\hR-\rohr}) --
                ({\d+\noTank*\dT+\noTank*\bT+\rV+\v+\v-\rohr},{\hR});
        \draw[]     ({\b},{\rohr}) --
                    ({\b},{\rohr*2+\hR+\hT}) --
                    ({\b+\dT+\rohr},{\rohr*2+\hR+\hT}) --
                    ({\b+\dT+\rohr},{\hR+\hT-\rohr*0.5}) ;
        \draw[]     ({\b+\rohr},{\rohr}) --
                    ({\b+\rohr},{\rohr+\hR+\hT}) --
                    ({\b+\dT},{\rohr+\hR+\hT}) --
                    ({\b+\dT},{\hR+\hT-\rohr*0.5}) ;
                    
        \filldraw[fill=white,draw=black] ({\b+\rohr*0.5},{\hR*2}) circle (\v*0.5);                    
        \fill[black]    ({\b+\rohr*1.2},{\hR*2}) --
                        ({\rohr},{\hR*2+\v*0.3}) --
                        ({\rohr},{\hR*2-\v*0.3}) --
                        ({\b+\rohr*1.2},{\hR*2}) 
                        coordinate (A);
                        

         \node at ({\dT +\d + 0.5*\bT},\hR+0.9*\hT) {\small tank 1};
         \node at ({2*\dT +\d + 1.5*\bT},\hR+0.9*\hT) {\small tank 2};
         \node at ({\bR*0.5},{\hR*0.7}) {reservoir};
         
         \draw[-latex] ({\b+0.5*\rohr},{\hR*2+\rohr}) --          ({\b+0.5*\rohr},\hR*3);
         
         \node at ({\d},{\hR*2.5}) {\small $u_1$} ;
         
         \draw[latex-latex] ({\d+\dT+0.5*\bT},{\hR}) --
                    ({\d+\dT+0.5*\bT},{\hR+0.7*\hT})
                    node [anchor = north east] {$x_1$};
        \draw[latex-latex] ({\d+2*\dT+1.5*\bT},{\hR}) --
                    ({\d+2*\dT+1.5*\bT},{\hR+0.5*\hT})
                    node [anchor = north east] {$x_2$};

        \draw[thin,dashed] ({\d+\noTank*\dT+1.5*\bT},{\hR+0.6*\hT}) ellipse ({0.5*\bT} and {\factor*0.1}) ;
        \node at  ({\d+\noTank*\dT+1.5*\bT},{\hR+0.7*\hT}) {\small $A$};
    \end{tikzpicture}
    \caption{Nonlinear water tank system. The tanks are connected by the valve $V_{12}$. Water can be pumped into the first tank with the control input $u_1(t)$ and is drained from the second tank with the valve $V_{2R}$.}
    \label{fig:watertank}
\end{figure}


\begin{table}[h!t]
    \begin{center}
    \caption{Parameters of the water tank system.}
        \begin{tabular}{llll}
            Parameter & Short form & Value & Unit \\ \hline
            Cross-sectional area &
            $A$ &
            $0.015$ &
            $\text{m}^2$ \\
            Maximum flow rate of pump one &
            $u_{1,{\text{max}}}$ &
            $2\cdot 10^{-4}$&
            ${\text{m}^3}/{\text{s}}$ \\
            Valve parameter $V_{12}$ &
            $c_{12}$ &
            $2.5\cdot 10^{-5}$&
            $\text{m}^2$ \\
            Valve parameter $V_{2R}$ &
            $c_{2\text{R}}$& 
            $2.5\cdot 10^{-5}$&
            $\text{m}^2$ \\ 
            Gravitational force &
            $g$ &
            $9.81$ &
            ${\text{m}}/{\text{s}^2}$ \\ 
            \hline	
        \end{tabular}
        \label{tab:parameters}
    \end{center}
\end{table}
To obtain a linearized state space representation, we first set $\dot{\state}_1(t), \dot{\state}_2(t)$ in Equation~\eqref{5_dreitankEq} equal to zero, solve for $\state_{e_1},\state_{e_2}$ depending on $\control_{e_1}$ and determine the Jacobian matrices $A_e(\state_e,\control_e)$, $B_e(\state_e,\control_e)$.
Now, we can choose a control input $\control_{e_1}$, calculate the associated equilibrium point $\state_e$ and insert it into $A_e(\state_e,\control_e)$, $B_e(\state_e,\control_e)$. For the given system we obtain unique solutions under the condition $\state(t) > 0$. However, this is not the case for every system and there may be multiple solutions or no solution at all.
From the linearized approximation of Equation~\eqref{5_dreitankEq} we obtain the Smith Normal Form
\begin{equation}
    D = \begin{bmatrix}
        1 & 0 & 0\\
        0 & 1 & 0
    \end{bmatrix},
\end{equation}
thus all differential operators of the system equations are contained in $V$.
Following the rules in \cite{besginow2022constraining}, we can construct the latent \ac{gp}
\begin{equation}
    h(t) \sim \mathcal{GP}\left(\begin{bmatrix}0 \\ 0 \\ 0 \end{bmatrix}, 
        \begin{bmatrix}
            0 & 0 & 0\\
            0 & 0 & 0\\
            0 & 0 & k_\text{SE}
        \end{bmatrix}\right)
\end{equation}
where $k_\text{SE}$ is the \ac{se} kernel, parameterized by $\sigma_f$ and $\ell$. 
The occurrence of one zero vector in $D$, and therefore one \ac{se} kernel in the covariance matrix, results from one degree of freedom in the system and is consistent with the single control input.

\subsection{Simulation Results}
According to the objective function Equation~\eqref{eq:mpc:obj} we investigate the mean control error
\begin{equation} \label{eq:control_error}
    \frac{1}{T}\sum_{i=1}^T (\state(t_i) - \state_{\text{ref}})^2
\end{equation}
and the mean control input
\begin{equation} \label{eq:control_sum}
    \frac{1}{T}\sum_{i=1}^T \Vert \control(t) \Vert
\end{equation}
in order to compare the control performance.
Furthermore, we investigate the mean constraint violation
\begin{equation}\label{eq:constr_viol}
    \frac{1}{T}\sum_{i=1}^T  \max\{z(t_i) - z_{\max}, 0\} + \max\{z_{\min} - z(t_i), 0\}
\end{equation}
to demonstrate, whether our approach can handle the imposed constraints.

Our control task is to track a constant reference at the equilibrium point $\state_{e,\text{ref}}$ with $\control_{e,\text{ref}}=0.3\cdot \control_{1,{\text{max}}}$, starting from the initial equilibrium point $\state_{e,0}$ with $\control_{e,0}=0.2\cdot \control_{1,{\text{max}}}$.
We compare three controller models (A), (B) and (C), which incorporate the reference in different ways.
We condition all models on the dataset $\mathcal{D}=\mathcal{D}_\text{init} \cup \mathcal{D}_\text{con}$, with $10$ equidistant datapoints from $t_1=1$s to $t_T=10$s for $\mathcal{D}_\text{con}$ according to Equation~\eqref{eq:soft_constraint}-\eqref{eq:soft_constraint_var}.
The hyperparameters are optimized offline in advance using $\mathcal{D}$ and the models are set up as follows:
\begin{itemize}
    \item [(A)]
    First, we set soft constraints  $\mathcal{D}_\text{con}$ according to their physical limit, that is 
    $\state(t) \in[0.0,0.6]^2$, $\control(t) \in[0.0, \control_{1,{\text{max}}}]$
    and include the reference as prior mean of the \ac{gp} $\mu_{\text{prior}}(t)$.
    \item [(B)]
    Next, we additionally incorporate the reference as endpoint constraint $\mathcal{D}_\text{ref}$ at time $t_{\text{ref}}=100$ s with zero noise.
    \item [(C)]
    At last, we place the soft constraints with 
    \[
    \state \in[0.9 \cdot \state_{e,1}, 1.1 \cdot \state_{e,1}]^\top, \ \control \in[0.9 \cdot \control_{e,1}, 1.1 \cdot \control_{e,1}]^\top
    \]
    close to the reference but don't add the endpoint constraint to the training data.
\end{itemize}

\begin{figure}[h]
    \centering
    \includegraphics[width=0.485\textwidth]{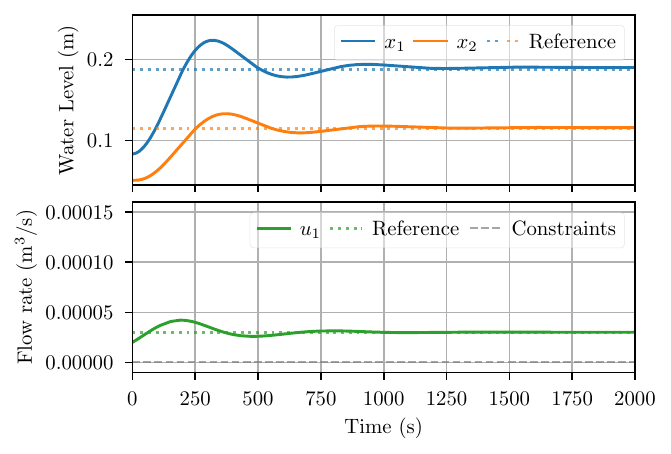}
    \vspace{-.8cm}
    \caption{Model (A): The reference is present in the \ac{gp} prior mean.
    Note the different time span.}
    \label{fig:result30}
\end{figure}

\begin{figure}[h]
    \centering
    \includegraphics[width=0.485\textwidth]{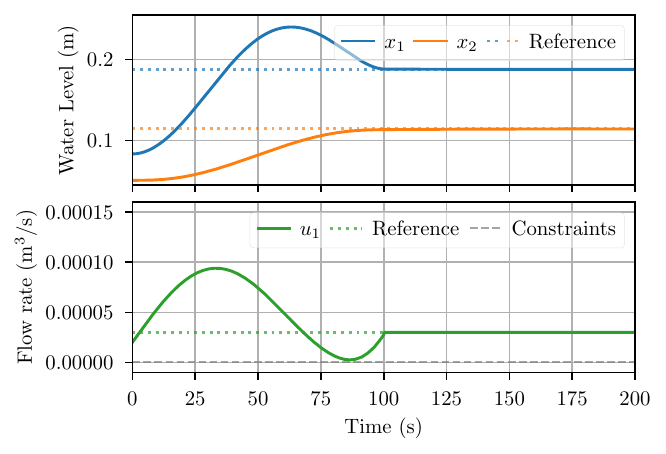}
    \vspace{-.8cm}
    \caption{Model (B): The reference is added to the training data as endpoint constraint at $t=100$ s.}
    \label{fig:result29}
\end{figure}

\begin{figure}[h]
    \centering
    \includegraphics[width=0.485\textwidth]{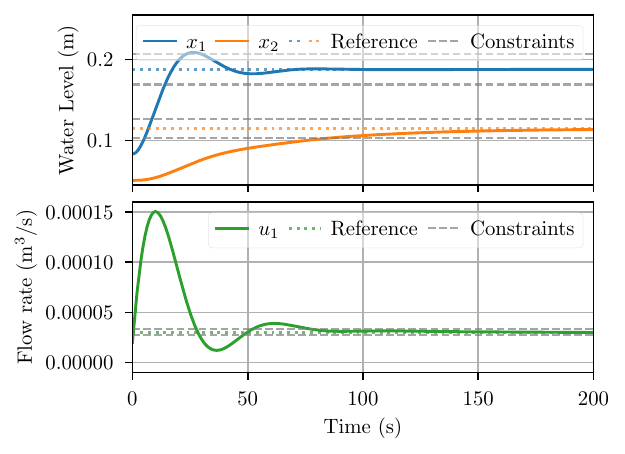}
    \vspace{-.8cm}
    \caption{Model (C): The reference is incorporated as soft constraints.}
    \label{fig:result27}
\end{figure}

Results are shown in Figure~\ref{fig:result30}-\ref{fig:result27} and Table~\ref{tab:regulation}.
Note the longer time span in Figure~\ref{fig:result30} in comparison to  Figure~\ref{fig:result29}-\ref{fig:result27}. 
Although the state trajectories converge to the reference in all models, the incorporation of additional information in the last two models (B) and (C) leads to noticeably faster convergence.

\begin{table}[ht]
    \centering
     \caption{Results for the control task for $200$ s.}
    \begin{tabular}{llll}
    Training dataset & Model (A) & Model (B) & Model (C)  \\
    \hline
    Control error       Eq.~\eqref{eq:control_error}
    & 5.05E-03	& 1.36E-03   & 4.19E-04 \\
    Mean Control input  Eq.~\eqref{eq:control_sum}
    & 3.40E-05  & 4.01E-05   & 1.18E-05 \\
    Constraint error    Eq.~\eqref{eq:constr_viol}
    & 0.0      & 0.0        & 4.3E-04 \\
    \hline
    \end{tabular}
    \label{tab:regulation}
\end{table}
The endpoint constraints in model (B) enforce the states to reach the reference at given time and after that, the states can be held constant, by applying $\control_{e,1}$ to the system. However, the soft constraints are almost violated in order to meet the hard constraints and enforcing a shorter convergence time would lead to control trajectories with negative values. 
While the higher overshoot of $\state_1(t)$ in model (B) achieves that all states converge at the same time, the controller in model (C) steers the trajectories not at the same time but smoother to the reference.
\section{Conclusion}
\label{sec:conclusion}
We implemented the novel \ac{lodegp}-based \ac{mpc} approach for nonlinear systems through linearization around an equilibrium point, which also served as a constant reference to achieve open-loop stability.
The reformulation of the control task as a \ac{gp} inference problem yields smooth control trajectories.
Hard and soft constraints can be incorporated as setpoints in the training data
by scaling the noise variance. 
To realize hard box constraints, we need to implement bounded likelihoods in the \ac{lodegp} as done by \cite{jensen_bounded_2013}.
At this point, we designed a local controller in the neighborhood of one equilibrium point. 
The combination of multiple local \acp{gp} is a well-researched field \citep{Nguyen-tuong_Peters_Seeger_2008, Gogolashvili_Kozyrskiy_Filippone_2022} 
and in the future we will extend our controller globally as it is done in gain scheduling \citep{Adamy_2022}.

\bibliography{ifacconf_2}

\begin{thebibliography}{37}
\providecommand{\natexlab}[1]{#1}
\providecommand{\url}[1]{\texttt{#1}}
\providecommand{\urlprefix}{URL }
\expandafter\ifx\csname urlstyle\endcsname\relax
  \providecommand{\doi}[1]{doi:\discretionary{}{}{}#1}\else
  \providecommand{\doi}{doi:\discretionary{}{}{}\begingroup
  \urlstyle{rm}\Url}\fi

\bibitem[{Adamy(2022)}]{Adamy_2022}
Adamy, J. (2022).
\newblock \emph{Nonlinear Systems and Controls}.
\newblock Springer Berlin Heidelberg, Berlin, Heidelberg.

\bibitem[{Berberich et~al.(2021)Berberich, Köhler, Müller, and
  Allgöwer}]{berberich_data-driven_2021}
Berberich, J., Köhler, J., Müller, M.A., and Allgöwer, F. (2021).
\newblock Data-driven model predictive control with stability and robustness
  guarantees.
\newblock \emph{{IEEE} Transactions on Automatic Control}, 66(4), 1702--1717.

\bibitem[{Besginow and Lange-Hegermann(2022)}]{besginow2022constraining}
Besginow, A. and Lange-Hegermann, M. (2022).
\newblock Constraining {G}aussian processes to systems of linear ordinary
  differential equations.
\newblock In \emph{NeurIPS}, volume~35.

\bibitem[{Chyzak et~al.(2005)Chyzak, Quadrat, and
  Robertz}]{chyzak2005effective}
Chyzak, F., Quadrat, A., and Robertz, D. (2005).
\newblock Effective algorithms for parametrizing linear control systems over
  ore algebras.
\newblock \emph{Applicable Algebra in Engineering, Communication and
  Computing}, 16, 319--376.

\bibitem[{Draeger et~al.(1995)Draeger, Engell, and Ranke}]{draeger_model_1995}
Draeger, A., Engell, S., and Ranke, H. (1995).
\newblock Model predictive control using neural networks.
\newblock \emph{{IEEE} Control Systems Magazine}, 15(5), 61--66.

\bibitem[{Gogolashvili et~al.(2022)Gogolashvili, Kozyrskiy, and
  Filippone}]{Gogolashvili_Kozyrskiy_Filippone_2022}
Gogolashvili, D., Kozyrskiy, B., and Filippone, M. (2022).
\newblock Locally smoothed gaussian process regression.
\newblock \emph{Procedia Computer Science}, 207, 2717–2726.

\bibitem[{Grüne and Pannek(2017)}]{Grüne_Pannek_2017}
Grüne, L. and Pannek, J. (2017).
\newblock \emph{Nonlinear Model Predictive Control}.
\newblock Communications and Control Engineering. Springer International
  Publishing, Cham.

\bibitem[{Harkonen et~al.(2023)Harkonen, Lange-Hegermann, and
  Raita}]{harkonen2023gaussian}
Harkonen, M., Lange-Hegermann, M., and Raita, B. (2023).
\newblock Gaussian process priors for systems of linear partial differential
  equations with constant coefficients.
\newblock In \emph{ICML}. PMLR.

\bibitem[{Hewing et~al.(2018)Hewing, Liniger, and
  Zeilinger}]{hewing2018cautious}
Hewing, L., Liniger, A., and Zeilinger, M.N. (2018).
\newblock Cautious nmpc with gaussian process dynamics for autonomous miniature
  race cars.
\newblock In \emph{ECC}, 1341--1348. IEEE.

\bibitem[{Jensen et~al.(2013)Jensen, Nielsen, and Larsen}]{jensen_bounded_2013}
Jensen, B.S., Nielsen, J.B., and Larsen, J. (2013).
\newblock Bounded gaussian process regression.
\newblock In \emph{International Workshop on Machine Learning for Signal
  Processing ({MLSP})}.

\bibitem[{Jidling et~al.(2017)Jidling, Wahlstr{\"o}m, Wills, and
  Sch{\"o}n}]{jidling2017linearly}
Jidling, C., Wahlstr{\"o}m, N., Wills, A., and Sch{\"o}n, T.B. (2017).
\newblock Linearly constrained gaussian processes.
\newblock \emph{NeurIPS}, 30.

\bibitem[{Kocijan et~al.(2004)Kocijan, Murray-Smith, Rasmussen, and
  Girard}]{Kocijan_Murray-Smith_Rasmussen_Girard_2004}
Kocijan, J., Murray-Smith, R., Rasmussen, C., and Girard, A. (2004).
\newblock Gaussian process model based predictive control.
\newblock In \emph{ACC}, 2214–2219 vol.3. IEEE, Boston, MA, USA.

\bibitem[{Kouvaritakis et~al.(2000)Kouvaritakis, Cannon, and
  Rossiter}]{Kouvaritakis_Cannon_Rossiter_2000}
Kouvaritakis, B., Cannon, M., and Rossiter, J.A. (2000).
\newblock Stability, feasibility, optimality and the degrees of freedom in
  constrained predictive control.
\newblock In \emph{Nonlinear Model Predictive Control}, 99–113. Birkhäuser
  Basel, Basel.

\bibitem[{Lange-Hegermann(2018)}]{langehegermann2018algorithmic}
Lange-Hegermann, M. (2018).
\newblock Algorithmic linearly constrained {G}aussian processes.
\newblock In \emph{NeurIPS}, volume~31.

\bibitem[{Lange-Hegermann and Robertz(2013)}]{lange2013thomas}
Lange-Hegermann, M. and Robertz, D. (2013).
\newblock Thomas decompositions of parametric nonlinear control systems.
\newblock \emph{IFAC Proceedings Volumes}.

\bibitem[{Lange-Hegermann and Robertz(2020)}]{lange2020thomas}
Lange-Hegermann, M. and Robertz, D. (2020).
\newblock Thomas decomposition and nonlinear control systems.
\newblock \emph{Algebraic and Symbolic Computation Methods in Dynamical
  Systems}.

\bibitem[{Levine(2018)}]{levine2018reinforcement}
Levine, S. (2018).
\newblock Reinforcement learning and control as probabilistic inference:
  Tutorial and review.

\bibitem[{Maiworm et~al.(2021)Maiworm, Limon, and
  Findeisen}]{maiworm2021online}
Maiworm, M., Limon, D., and Findeisen, R. (2021).
\newblock Online learning-based model predictive control with gaussian process
  models and stability guarantees.
\newblock \emph{International Journal of Robust and Nonlinear Control}, 31(18),
  8785--8812.

\bibitem[{Matschek et~al.(2020)Matschek, Himmel, Sundmacher, and
  Findeisen}]{Matschek_Himmel_Sundmacher_Findeisen_2020}
Matschek, J., Himmel, A., Sundmacher, K., and Findeisen, R. (2020).
\newblock Constrained gaussian process learning for model predictive control.
\newblock \emph{IFAC-PapersOnLine}, 53(2), 971–976.

\bibitem[{Newman(1997)}]{newman1997smith}
Newman, M. (1997).
\newblock The {S}mith normal form.
\newblock \emph{Linear algebra and its applications}, 254(1-3), 367--381.

\bibitem[{Nguyen-Tuong et~al.(2008)Nguyen-Tuong, Peters, and
  Seeger}]{Nguyen-tuong_Peters_Seeger_2008}
Nguyen-Tuong, D., Peters, J., and Seeger, M. (2008).
\newblock Local gaussian process regression for real time online model
  learning.
\newblock In \emph{NeurIPS}, volume~21. Curran Associates, Inc.

\bibitem[{Oberst(1990)}]{oberst1990multidimensional}
Oberst, U. (1990).
\newblock Multidimensional constant linear systems.
\newblock \emph{Acta Applicandae Mathematica}, 20(1), 1--175.

\bibitem[{Piche et~al.(2000)Piche, Sayyar-Rodsari, Johnson, and
  Gerules}]{piche_nonlinear_2000}
Piche, S., Sayyar-Rodsari, B., Johnson, D., and Gerules, M. (2000).
\newblock Nonlinear model predictive control using neural networks.
\newblock \emph{{IEEE} Control Systems Magazine}, 20(3), 53--62.

\bibitem[{Pommaret and Quadrat(1999)}]{pommaret1999algebraic}
Pommaret, J. and Quadrat, A. (1999).
\newblock Algebraic analysis of linear multidimensional control systems.
\newblock \emph{IMA Journal of Mathematical control and Information}, 16(3),
  275--297.

\bibitem[{Qin and Badgwell(2000)}]{Qin_Badgwell_2000}
Qin, S.J. and Badgwell, T.A. (2000).
\newblock An overview of nonlinear model predictive control applications.
\newblock In \emph{Nonlinear Model Predictive Control}, 369–392. Birkhäuser
  Basel.

\bibitem[{Rasmussen et~al.(2006)Rasmussen, Williams
  et~al.}]{rasmussen2006gaussian}
Rasmussen, C.E., Williams, C.K., et~al. (2006).
\newblock \emph{Gaussian processes for machine learning}, volume~1.
\newblock Springer.

\bibitem[{Rawlings et~al.(2017)Rawlings, Mayne, and Diehl}]{rawlings2017model}
Rawlings, J.B., Mayne, D.Q., and Diehl, M. (2017).
\newblock \emph{Model predictive control: theory, computation, and design},
  volume~2.
\newblock Nob Hill Publishing Madison, WI.

\bibitem[{Ross et~al.(2021)Ross, Smith, and Álvarez}]{Ross_Smith_Alvarez_2021}
Ross, M., Smith, M.T., and Álvarez, M. (2021).
\newblock Learning nonparametric volterra kernels with gaussian processes.
\newblock In \emph{NeurIPS}, volume~34, 24099--24110.

\bibitem[{Smith(1862)}]{smith1862systems}
Smith, H.J.S. (1862).
\newblock I. on systems of linear indeterminate equations and congruences.
\newblock \emph{Proceedings of the Royal Society of London}, (11), 86--89.

\bibitem[{Tebbe et~al.(2023)Tebbe, Pawlik, Trilling, L{\"o}bner,
  Lange-Hegermann, and Schneider}]{tebbe2023holistic}
Tebbe, J., Pawlik, T., Trilling, M., L{\"o}bner, J., Lange-Hegermann, M., and
  Schneider, J. (2023).
\newblock Holistic optimization of a dynamic cross-flow filtration process
  towards a cyber-physical system.
\newblock In \emph{INDIN}. IEEE.

\bibitem[{Tebbe et~al.(2025)Tebbe, Besginow, and
  Lange-Hegermann}]{tebbe2024physicsinformedgaussianprocesseslinear}
Tebbe, J., Besginow, A., and Lange-Hegermann, M. (2025).
\newblock Physics-informed gaussian processes as linear model predictive
  controller.
\newblock In \emph{L4DC}. PMLR.

\bibitem[{Torrisi et~al.(2016)Torrisi, Grammatico, Smith, and
  Morari}]{torrisi_projected_2016}
Torrisi, G., Grammatico, S., Smith, R.S., and Morari, M. (2016).
\newblock A projected gradient and constraint linearization method for
  nonlinear model predictive control.

\bibitem[{Umlauft et~al.(2018)Umlauft, Beckers, and
  Hirche}]{Umlauft_Beckers_Hirche_2018}
Umlauft, J., Beckers, T., and Hirche, S. (2018).
\newblock Scenario-based optimal control for gaussian process state space
  models.
\newblock In \emph{ECC}, 1386–1392. IEEE.

\bibitem[{Willems and Polderman(1997)}]{willems1997introduction}
Willems, J.C. and Polderman, J.W. (1997).
\newblock \emph{Introduction to mathematical systems theory: a behavioral
  approach}, volume~26.
\newblock Springer Science \& Business Media.

\bibitem[{Zerz(2000)}]{zerz2000topics}
Zerz, E. (2000).
\newblock \emph{Topics in multidimensional linear systems theory}, volume 256.
\newblock Springer Science \& Business Media.

\bibitem[{Zheng(2000)}]{Zheng_2000}
Zheng, A. (2000).
\newblock Some practical issues and possible solutions for nonlinear model
  predictive control.
\newblock In \emph{Nonlinear Model Predictive Control}, 129–143. Basel.

\bibitem[{Álvarez et~al.(2009)Álvarez, Luengo, and
  Lawrence}]{pmlr-v5-alvarez09a}
Álvarez, M., Luengo, D., and Lawrence, N.D. (2009).
\newblock Latent force models.
\newblock In \emph{AISTATS}, 9--16. PMLR.

\end{thebibliography}

\end{document}